\documentclass[11pt, a4paper, leqno]{amsart}

\usepackage{amsmath,amsthm,amssymb} 
\usepackage{amsfonts} 
\usepackage{mathrsfs} 
\usepackage[T1]{fontenc} 
\usepackage[utf8]{inputenc} 

\usepackage[dvipsnames, svgnames]{xcolor} 
\usepackage{enumerate} 
\usepackage[colorlinks,citecolor=citec,hypertexnames=false,urlcolor=urlc,breaklinks=true, pagebackref]{hyperref} 
\usepackage{thmtools} 
\usepackage{accents} 

\usepackage{amsbsy} 
\usepackage{amscd} 
\usepackage{latexsym} 
\usepackage{txfonts} 
\usepackage{exscale} 
\usepackage{bbm} 
 
\usepackage{geometry} 

\usepackage{cite} 

\newif\ifsoldark
\newif\ifsollight
\newif\ifclassic
\newif\ifplain

\numberwithin{equation}{section}

\declaretheoremstyle[
headfont=\color{lmcolor}\normalfont\bfseries,
bodyfont=\normalfont\itshape,
]{colorlemma}
\theoremstyle{colorlemma}

\declaretheoremstyle[
headfont=\color{propcolor}\normalfont\bfseries,
bodyfont=\normalfont\itshape,
]{colorprop}
\theoremstyle{colorprop}
\newtheorem{proposition}[equation]{Proposition}

\declaretheoremstyle[
headfont=\color{thmcolor}\normalfont\bfseries,
bodyfont=\normalfont\itshape,
]{colorthm}
\theoremstyle{colorthm}
\newtheorem{theorem}[equation]{Theorem}

\declaretheoremstyle[
headfont=\color{defcolor}\normalfont\bfseries,
bodyfont=\normalfont,
]{colordef}
\theoremstyle{colordef}

\declaretheoremstyle[
headfont=\color{excolor}\normalfont\bfseries,
bodyfont=\normalfont,
]{colorexample}
\theoremstyle{colorexample}

\declaretheoremstyle[
headfont=\color{rmcolor}\normalfont\itshape,
bodyfont=\normalfont,
]{colorremark}
\theoremstyle{colorremark}
\newtheorem{remark}[equation]{Remark}
\newtheorem*{remark*}{Remark}

\makeatletter
\renewcommand{\eqref}[1]{\textup{\eqreftagform@{\ref{#1}}}}
\let\eqreftagform@\tagform@
\def\tagform@#1{%
	\maketag@@@{\color{eqcolor}(\ignorespaces#1\unskip\@@italiccorr)}%
}
\makeatother

\newcommand{\dint}{\int\!\!\!\int}

\newcommand{\dv}{\operatorname{div}}
\newcommand{\n}[1]{\mathscr{#1}}

\newcommand{\bb}[1]{\mathbb{#1}}

\newcommand{\RNum}[1]{\uppercase\expandafter{\romannumeral #1\relax}}

\DeclareMathOperator{\diam}{diam}

\def\XXint#1#2#3{{\setbox0=\hbox{$#1{#2#3}{\int}$}
		\vcenter{\hbox{$#2#3$}}\kern-.5\wd0}}

\def\YYint#1#2#3{{\setbox0=\hbox{$#1{#2#3}{\iint}$}
		\vcenter{\hbox{$#2#3$}}\kern-.51\wd0}}
 
\newcommand{\ra}{\rightarrow}

\newcommand{\m}[1]{\mathcal{#1}}

\newcommand{\vertiii}[1]{{\left\vert\kern-0.25ex\left\vert\kern-0.25ex\left\vert #1 
		\right\vert\kern-0.25ex\right\vert\kern-0.25ex\right\vert}}

\allowdisplaybreaks

\calclayout 

\classictrue 

\ifsoldark

\definecolor{solblack}{HTML}{001116} 	
\definecolor{solnormal}{HTML}{b3cacc}	
\definecolor{solblue}{HTML}{268bd2}
\definecolor{solpurple}{HTML}{6c71c4}
\definecolor{solteal}{HTML}{2aa198}
\definecolor{solgreen}{HTML}{859900}
\definecolor{solpink}{HTML}{ffaaff}
\definecolor{solbrown}{HTML}{b58900}
\definecolor{solorange}{HTML}{cb4b16}
\definecolor{solviolet}{HTML}{8d95ff}

\colorlet{citec}{solblue}
\colorlet{urlc}{solteal}
\colorlet{toc}{solpink}
\colorlet{hyperc}{solviolet}
\colorlet{bpcolor}{solblue}
\colorlet{smcolor}{solpurple}
\colorlet{impcolor}{solpink}
\colorlet{eqcolor}{solgreen}
\colorlet{lmcolor}{solbrown}
\colorlet{propcolor}{solpurple}
\colorlet{thmcolor}{solviolet}
\colorlet{defcolor}{solblue}
\colorlet{rmcolor}{solorange}
\colorlet{excolor}{solgreen}
\pagecolor{solblack}
{\color{solnormal} 
	
	\fi
	
	\ifsollight
	
	solarized colors
	\definecolor{solblack}{HTML}{002b36}
	\definecolor{solnormal}{HTML}{839496}
	\definecolor{solblue}{HTML}{268bd2}
	\definecolor{solpurple}{HTML}{6c71c4}
	\definecolor{solteal}{HTML}{2aa198}
	\definecolor{solgreen}{HTML}{859900}
	\definecolor{solpink}{HTML}{ffaaff}
	\definecolor{solbrown}{HTML}{b58900}
	\definecolor{solorange}{HTML}{cb4b16}
	\definecolor{solviolet}{HTML}{8d95ff}
	
	\fi
	
	\ifclassic
	
	\colorlet{citec}{blue}
	\colorlet{urlc}{teal}
	\colorlet{toc}{NavyBlue}
	\colorlet{hyperc}{Green}
	\colorlet{bpcolor}{NavyBlue}
	\colorlet{smcolor}{purple}
	\colorlet{impcolor}{ProcessBlue}
	\colorlet{eqcolor}{ForestGreen}
	\colorlet{lmcolor}{Mahogany}
	\colorlet{propcolor}{RoyalPurple}
	\colorlet{thmcolor}{NavyBlue}
	\colorlet{defcolor}{Cerulean}
	\colorlet{rmcolor}{MidnightBlue}
	\colorlet{excolor}{PineGreen}
	
	\fi

	\ifplain
	
	\colorlet{citec}{blue}
	\colorlet{urlc}{blue}
	\colorlet{toc}{red}
	\colorlet{hyperc}{blue}
	\colorlet{bpcolor}{NavyBlue}
	\colorlet{smcolor}{purple}
	\colorlet{impcolor}{ProcessBlue}
	\colorlet{eqcolor}{black}
	\colorlet{lmcolor}{black}
	\colorlet{propcolor}{black}
	\colorlet{thmcolor}{black}
	\colorlet{defcolor}{black}
	\colorlet{rmcolor}{black}
	\colorlet{excolor}{black}
	
	\fi

\begin{document}
	
\title[Kenig-Pipher condition and absolute continuity]{Failure to slide: a brief note on the interplay between the Kenig-Pipher condition and the absolute continuity of elliptic measures}

\author{B. Poggi}

\address{Bruno Giuseppe Poggi Cevallos
	\\
	School of Mathematics
	\\
	University of Minnesota
	\\
	Minneapolis, MN 55455, USA} \email{poggi008@umn.edu}

\maketitle

\date{\today}

\keywords{}

\begin{abstract}
In this   note, we explore some consequences of the Modica-Mortola construction of a singular elliptic measure, as regards the link between the quantitative absolute continuity ($A_{\infty}$) of their approximations and the suitability of a well-known tool, the so-called Kenig-Pipher condition ($\operatorname{KP}$). The Kenig-Pipher condition is used to ascertain absolute continuity in the presence of some mild regularity of the coefficient matrix. We perform some modifications of the Modica-Mortola example to show the following two statements: (a) there are sequences of matrices for which both $\operatorname{KP}$ and the $A_{\infty}$ condition break  down in the limit. (b) there are sequences of matrices for which $\operatorname{KP}$ breaks down but $A_{\infty}$ is preserved in the limit.
\end{abstract}

{
	\hypersetup{linkcolor=toc}
	\tableofcontents
}
\hypersetup{linkcolor=hyperc}

\section{Introduction} 

In this note, we explore some consequences of the Modica-Mortola construction \cite{mm} of a singular elliptic measure, as regards the link between the quantitative absolute continuity of their approximations and the suitability of a well-known tool, the so-called Kenig-Pipher condition.  

The Modica-Mortola example, in conjunction with the Caffarelli-Fabes-Kenig example \cite{cfk}, were the first (and concurrent) constructions of an elliptic measure singular with respect to the surface measure of a smooth domain. The former uses an approximation procedure, lacunary sequences, and Riesz products, while the latter relies on the theory of quasi-conformal mappings. The interest in evidencing such cases had been aroused since Dahlberg \cite{dah1} proved a few years earlier that the elliptic measure for the Laplacian $L=-\Delta$ was absolutely continuous with respect to the surface measure of the unit ball. On the other hand, Caffarelli-Fabes-Mortola-Salsa \cite{cfms} had shown that all elliptic measures were doubling, precluding the existence of trivial examples.

Ever since then, understanding the precise relationship between the coefficients of a divergence-form elliptic operator and the absolute continuity of the elliptic measure has been an ongoing and lively area of research, whose review we defer to any one of the many contemporaneous papers in the landscape. We do bring attention to one of the landmark results in the literature, \cite{kp3}, in which it is shown that quantifiable absolute continuity of the elliptic measure on the unit ball can be ascertained when the gradient of the coefficient matrix satisfies a Carleson-measure type condition. Their condition has come to be known as the \emph{Kenig-Pipher condition}, and its connection to the absolute continuity of elliptic measures has been seen to be remarkably robust in the past decade.

Our aim here is modest. We will adapt the game played by Modica-Mortola to
\begin{enumerate}[(a)]
	\item provide a sequence of $A_{\infty}$ elliptic measures for whom their Kenig-Pipher condition breaks down in the limit, but which nevertheless converge weakly to a singular elliptic measure, and
	\item provide a sequence of $A_{\infty}$ elliptic measures for whom their Kenig-Pipher condition breaks down in the limit, but which converge to an absolutely continuous elliptic measure.
\end{enumerate}
For case (a), we also argue that the  placement in $A_{\infty}$ of the approximating measures degenerates. This shows that, in a sense (see Remark \ref{rm.top}), the Modica-Mortola matrix   lies at the boundary of the set of matrices with $A_{\infty}$  elliptic measures. 

\section{The Modica-Mortola example: an exposition}\label{sec.mm}

Let us first supply a brief exposition of the construction in \cite{mm}. In preparation, note that if $A, B$ are two positive functions, we say that $A\approx B$ if there exists a constant $C\geq1$ such that
\[
\frac1CB\leq A\leq C B.
\]
If $C$ depends on some parameter $\beta$, then we make the dependence on $\beta$ explicit by using $A\approx_{\beta} B$ instead. We reserve the notation $dx$ for the Lebesgue measure on an interval.

The idea of the construction is to jam the boundary with very thin layers consisting of a material with highly oscillatory periodic anisotropy.

Let
\[
A(x,y)=\begin{pmatrix}1&0\\0&\alpha(x,y)\end{pmatrix},\qquad A_j(x,y)=\begin{pmatrix}1&0\\0&\alpha_j(x,y)\end{pmatrix},~j\in\bb N,~ (x,y)\in\bb R^2,
\]
where
\[
\alpha(x,y)=\begin{cases}\phi_1(x),\qquad\text{if }|y|\geq1/k_1,\\[2mm]\psi(k_{j+1}y)\phi_{j+1}(x)+(1-\psi(k_{j+1}y))\phi_j(x),\qquad\text{if }\frac1{k_{j+1}}\leq|y|<\frac1{k_j}, j=1,2,\ldots,\\[2mm]1,\qquad\text{if }y=0,\end{cases}
\]
\[
\alpha_j(x,y)=\begin{cases}\alpha(x,y),\qquad\text{if }|y|\geq1/k_j,\\[2mm]\phi_j(x)\qquad\text{if }|y|<1/k_j,\end{cases}
\]
\begin{equation}\label{eq.phij}
\phi_j(x)=1+\frac1{4\pi\sqrt j}\cos(2\pi h_jx),\qquad j\in\bb N,
\end{equation}
and $\{h_j\}$, $\{k_j\}$ are suitable lacunary sequences (see, for instance, \cite{zyg}) of positive integers which can be chosen to satisfy that 
\begin{equation}\label{eq.lac}
h_{j+1}\geq4h_j,\qquad k_{j+1}\geq2k_j, \qquad\text{and }  h_j\geq jk_j\footnote{The condition $h_j\geq jk_j$ does not appear in \cite{mm}; however, it is clear by their inductive construction of the sequences in pages 16-17 that we may choose the sequences in this way, since we choose $h_{j+1}$ large based on already having chosen $k_{j+1}$. Note that asking this condition simplifies the proof of Proposition \ref{prop.break} below, but is probably not strictly required to achieve the result.},
\end{equation}
and $\psi\in C_c^{\infty}(\bb R)$ is a smooth cut-off satisfying $\psi(t)=\psi(-t), 0\leq\psi(t)\leq1$ for all $t\in\bb R$, $\psi(t)=1$ if $|t|\leq1$, and $\psi(t)=0$ if $|t|\geq2$.

Observe that $\alpha\in C^0(\bb R^2)\cap C^{\infty}(\{(x,y)\in\bb R^2:y\neq0\})$ and $\frac12\leq\alpha(x,y)\leq\frac32$ for every $(x,y)\in\bb R^2$. Note that $\alpha_j\in C^{\infty}(\bb R^2)$ for each $j\in\bb N$ and $\alpha_j$ converges pointwise uniformly in $\bb R^2$ to $\alpha$.

Let $\Omega$ be an open, bounded subset of the upper-half plane $\bb R\times\bb R_+$, with smooth boundary and such that
\[
[-10,20]\times\{0\}\subset\partial\Omega,
\]
and let $L_j:=-\operatorname{div }A_j\nabla$ and $L:=-\operatorname{div }A\nabla$ be operators on $\Omega$ with corresponding elliptic measures $\{\omega_j^P\}_{P\in\Omega}$ and $\{\omega^P\}_{P\in\Omega}$ on $\partial\Omega$. Henceforth, we fix an arbitrary $P\in\Omega$ and write $\omega_j=\omega_j^P$, $\omega=\omega^P$. We say that a measure $\mu$ is \emph{singular} with respect to a measure $\nu$ if $\mu$ is not absolutely continuous with respect to $\nu$. The main result in \cite{mm} is

\begin{theorem}[A singular elliptic measure; \cite{mm}]\label{thm.main} The probability measure $\omega$ on $\partial\Omega$ is singular with respect to the surface measure.
\end{theorem}

\noindent\emph{Proof.}  In this situation, $\omega_j^P$ converges weakly in the sense of measures to $\omega^P$ on $\partial\Omega$ as $j\ra\infty$ (cf. Lemma 1 in \cite{mm}).  For each $j\in\bb N$, let $g_j=g_j(\cdot,P)$ be the \emph{Green function} for the Dirichlet problem for the operator $L_j$ on $\Omega$, with pole $P$.  By the Green representation formula, for each $\chi\in C(\partial\Omega)$  we have the identity
\begin{equation}\label{eq.greenformula2}
\int_{\partial\Omega}\chi\,d\omega_j=\int_{\partial\Omega}\chi (A_j\nabla g_j)\cdot\hat n\,d\sigma,
\end{equation}
where $\hat n$ is the unit normal vector on $\partial\Omega$ pointing inward and $\sigma$ is the surface measure on $\partial\Omega$, which is the $(n-1)$-dimensional Hausdorff measure on $\partial\Omega$ (up to a dimensional constant).

Since $\Omega$ is a smooth bounded domain and each $A_j$ is smooth in $\bb R^2$, the classical Hopf lemma and the fact that $\alpha_j\approx 1$ imply that $(A\nabla g_j)\cdot\hat n\approx_j1$, so that $\omega_j$ and $\sigma$ are mutually absolutely continuous. Let $\m K_j=\m K_j^P$ be the \emph{Poisson kernel} of $\omega_j$ with pole $P$, which is the Radon-Nikodym derivative of $\omega_j$ with respect to $\sigma$. We have concluded that
\[
\m K_j\approx_j1,\qquad\text{on }\partial\Omega.
\]

Observe that by specializing (\ref{eq.greenformula2}) to the case when $\chi\in C_c([-1,1]\times\{0\})$, we procure the identity
\begin{equation}\label{eq.greenformula1}
\int_{\partial\Omega}\chi\,d\omega_j=\int_{-1}^1\chi(x,0)\phi_j(x)\frac{\partial g_j}{\partial y}(x,0)\,dx,
\end{equation}
so that pointwise almost everywhere on $[-1,1]\times\{0\}$, we have the representation
\[
\m K_j=\phi_j\frac{\partial g_j}{\partial y}.
\]
We will now sketch the fact that for each $j\geq2$, we may choose the lacunary sequences $\{h_j\}$, $\{k_j\}$ so that
\begin{equation}\label{eq.likeriesz}
\frac{\partial g_j}{\partial y}(x,0)\approx\prod_{i=1}^{j-1}\phi_i(x)=:\m R_{j-1}(x),\qquad x\in[-1,1],
\end{equation}
whence we summarily deduce that
\begin{equation}\label{eq.likerieszk}
\m K_j\approx\m R_j
\end{equation}
on $[-1,1]\times\{0\}$. We call $\m R_j$ a \emph{Riesz product} (see \cite{zyg}, Chapter V, Section 7). Let
\[
M:=\min_{|x|\leq1}\frac{\partial g_1}{\partial y}(x,0)
\]
and observe that, by the Hopf lemma, $M>0$. The sequences $\{h_j\}, \{k_j\}$ are chosen so that the estimate
\begin{equation}\label{eq.greendecay}
\max_{|x|\leq1}\Big|\frac{\partial g_{j+1}}{\partial y}(x,0)-\phi_j(x)\frac{\partial g_j}{\partial y}(x,0)\Big|\leq M4^{-j-1}
\end{equation}
holds for all $j\in\bb N$. That this can be done is much of the program in \cite{mm}, and thus we leave the study of this technology to them. Consider for each $j\in\bb N$ the function
\[
w_j(x):=\frac{\frac{\partial g_{j+1}}{\partial y}(x,0)}{\m R_j(x)},\qquad |x|\leq 1,~j\in\bb N,
\]
which can easily be rewritten as
\[
w_j(x)=\frac{\partial g_1}{\partial y}(x,0)+\sum_{i=1}^j\frac{\frac{\partial g_{i+1}}{\partial y}(x,0)-\phi_i(x)\frac{\partial g_i}{\partial y}(x,0)}{\m R_i(x)}.
\]
From the above equality, (\ref{eq.greendecay}), and the fact that $\phi_j\geq\frac12$ for each $j\in\bb N$, it follows that $\{w_j\}$ is a Cauchy sequence in $C[-1,1]$, whence there exists $w\in C[-1,1]$ so that $w_j\ra w$ in $C[-1,1]$, and moreover
\[
\max_{|x|\leq1}|w_j(x)-w(x)|\leq M2^{-n-2},\qquad\text{for each }j\in\bb N, 
\]
\[
w\geq\frac34M,\qquad\text{on }[-1,1].
\]
These computations prove (\ref{eq.likeriesz}) and therefore (\ref{eq.likerieszk}). We now borrow from \cite{zyg} Chapter V, Section 7, Lemma 7.5, the wisdom that these \emph{Riesz products are the partial sums of a Fourier-Stieltjes series of a non-decreasing (non-constant) continuous function $F$ on $[-1,1]$, whose derivative is $0$ almost everywhere}. In particular, $\m R_j\,dx$ converges weakly in the sense of  measures to a singular measure $dF$ on $[-1,1]$. Since we have (\ref{eq.likerieszk}) and the fact that $\omega_j\ra\omega$ weakly in the sense of measures, it follows that $\omega$ is singular, as desired.\hfill{$\square$}

\section{The Modica-Mortola approximations, their absolute continuity, and the Kenig-Pipher condition}

Let $\Omega\subset\bb R^2$ be as above. Given $X\in\Omega$, denote by $\delta(X)$ the distance from $X$ to $\partial\Omega$. For an elliptic matrix $\m A$ of bounded real measurable coefficients in $\Omega$, we say that $\m A$ satisfies the \emph{Kenig-Pipher condition} if the quantity
\[
\n P(\m A):=\sup_{{\tiny\begin{matrix}q\in\partial\Omega\\0<r<\operatorname{diam}(\Omega)\end{matrix}}}\frac1r\dint_{B(q,r)\cap\Omega}\Big(\sup_{Y\in B(X,\frac{\delta(X)}2)}|\nabla\m A(Y)|^2\delta(Y)\Big)\,dX 
\]
is finite. It turns out that if the matrix $\m A$ satisfies the Kenig-Pipher condition, then the elliptic measure $\omega_{\m A}$ associated to the operator $L=-\dv\m A\nabla$ is absolutely continuous with respect to the surface measure. In fact, one can say more: in this case, the absolute continuity can be \emph{quantified} using the theory of Muckenhoupt $A_p$ weights (see \cite{ste} for definitions and the basic results). 

Thus, it can be shown that if $\m A$ satisfies the Kenig-Pipher condition, then $\omega_{\m A}\in A_{\infty}(\sigma)$. We now briefly summarize some results regarding the Muckenhoupt $A_{\infty}$ class and Reverse H\"older classes which we shall later use. The characterization $\omega_{\m A}\in A_{\infty}(\sigma)$ is equivalent to (and, therefore, may be defined as) the condition  that $\omega_{\m A}$ is absolutely continuous with respect to the surface measure $\sigma$ and such  that the Poisson kernel $\m K_{\m A}$ lies in $RH_q$ for some $q>1$, where $RH_q$ is the space of non-negative weights satisfying a \emph{Reverse-H\"older inequality}: a non-negative weight $w$ on $\partial\Omega$  lies in $RH_q$ for $q>1$ if there exists a constant $C\geq1$ such that for all surface balls $\Delta=B\cap\partial\Omega$ ($B$ is an $n-$dimensional ball centered on $\partial\Omega$), the estimate
\begin{equation}\label{eq.rhcond}
\Big(\frac1{\sigma(\Delta)}\int_{\Delta}w^q\,d\sigma\Big)^{\frac1q}\leq C\frac1{\sigma(\Delta)}\int_{\Delta}w\,d\sigma 
\end{equation}
holds. We let $|w|_{RH_q}$ be the infimum of the set of all possible constants $C$ such that (\ref{eq.rhcond}) holds. Moreover, we have the following direct characterization of $A_{\infty}$ weights\footnote{here, there is a slight abuse of notation, as we consider both $A_{\infty}$ for measures and for weights; these are essentially equivalent, however.}: a non-negative weight $w$ lies in $A_{\infty}$ if and only if the quantity
\begin{equation}
|w|_{A_{\infty}}:=\sup_{\Delta\subset\partial\Omega}\Big\{\Big(\frac1{\sigma(\Delta)}\int_{\Delta}w\,d\sigma\Big)\exp\Big(\frac1{\sigma(\Delta)}\int_{\Delta}\log w^{-1}\,d\sigma\Big)\Big\}
\end{equation}
is finite. We call $|w|_{A_\infty}$ the \emph{$A_{\infty}$ constant of w}.  The limiting case of the Reverse H\"older classes is the space $RH_{L\log L}$ of weights which satisfy the reverse Jensen's inequality for the function $x\log x$:
\[
\Vert w\Vert_{(L\log L,\frac{d\sigma|_{\Delta}}{\sigma(\Delta)})}\leq C\frac1{\sigma(\Delta)}\int_{\Delta}w\,d\sigma,
\]
and we call $|w|_{RH_{L\log L}}$ the infimum of the set of all possible $C$ such that the above inequality holds. We have that $A_{\infty}=RH_{L\log L}=\bigcup_{q>1}RH_q$.

Using the notation of the previous section, note that for each fixed $j\in\bb N$, $A_j$ trivially satisfies the Kenig-Pipher condition. Indeed, since $A_j$ is smooth on $\bb R^2$, there exists a constant $C_j$ such that $|\nabla A_j|\leq C_j$ on $\overline{\Omega}$, whence 
\[
\n P(A_j)\lesssim_j\diam(\Omega)^2<\infty.
\]
This fact provides a second easy proof of the fact that for each $j\in\bb N$, $\omega_j\in A_{\infty}(\sigma)$. 

The main calculations of this note follow. First, we check directly that the approximations $A_j$ of the Modica-Mortola example break the Kenig-Pipher condition ``in the limit''.


\begin{proposition}[Breaking of the Kenig-Pipher condition on a sliding scale]\label{prop.break} There exists a sequence $\{A_j\}$ of diagonal elliptic matrices on $\Omega$ with smooth, bounded, real coefficients in $\Omega$ and uniformly continuous on $\overline\Omega$, and there exists a diagonal elliptic matrix $A$ on $\Omega$ with smooth, bounded, real coefficients in $\Omega$ and uniformly continuous on $\overline\Omega$, such that $A_j\ra A$ pointwise uniformly on $\overline\Omega$, $\omega_A$ is singular with respect to the surface measure, and
\begin{equation}\label{eq.result}
\n P(A_j)\gtrsim j.
\end{equation}
\end{proposition}
\noindent\emph{Proof.} Let $\{A_j\}$, $A$ and all other variables be defined as in the previous section. Without loss of generality we may assume that $\Omega$ contains the square $[-1,2]\times[0,3]$. Fix $j\in\bb N$, and reckon the elementary estimates
\begin{multline}\nonumber
\n P(A_j)\geq\sup_{{\tiny\begin{matrix}x_0\in[0,1]\\0<r\leq\frac1{2k_j}\end{matrix}}}\frac1r\int_{x_0-r}^{x_0+r}\int_0^{\sqrt{r^2-(x-x_0)^2}}\Big(\sup_{(x_2,y_2)\in B((x,y),\frac y2)}\Big|\frac{\partial\alpha_j}{\partial x}(x_2,y_2)\Big|^2y_2\Big)\,dy\,dx\\ \geq\frac{h_j^2}{8j} \sup_{x_0\in[0,1]}\frac1{r_j}\int_{x_0-r_j}^{x_0+r_j}\int_0^{\sqrt{r_j^2-(x-x_0)^2}}y\Big(\sup_{x_2\in(x-\frac y2,\, x+\frac y2)}|\sin(h_jx_2)|^2 \Big)\,dy\,dx,
\end{multline}
where $r_j:=\frac1{2k_j}$. We may assume that both $h_j$ and $k_j$ are large, and $h_j\gg k_j$ (a fact afforded by virtue of the choice $h_j\geq jk_j$). Given $x_0\in[0,1]$, observe that
\[
\sup_{x_2\in(x-\frac y2,\, x+\frac y2)}|\sin(h_jx_2)|^2 =1
\]
for each $(x,y)$ in the set
\begin{multline*}
S:=\Big\{(x,y):x\in(x_0-r_j,x_0-\tfrac{r_j}2)~,~y\in\big(\tfrac12\sqrt{r_j^2-(x-x_0)^2},\sqrt{r_j^2-(x-x_0)^2}\big)\Big\}\\ \subset\big\{(x,y):x\in(x_0-r_j,x_0+r_j)~,~y\in\big(0,\sqrt{r_j^2-(x-x_0)^2}\big)\big\},
\end{multline*}
because $\sin(h_j\cdot)$ is oscillating rapidly in an interval of length roughly $1/k_j$. Therefore,
\begin{multline}
\n P(A_j)\geq\frac{h_j^2}{8jr_j} \sup_{x_0\in[0,1]} \int_{x_0-r_j}^{x_0-\frac12r_j}\int_{\frac12\sqrt{r_j^2-(x-x_0)^2}}^{\sqrt{r_j^2-(x-x_0)^2}}\,y\,dy\,dx\\ \geq\frac{3h_j^2}{64jr_j} \sup_{x_0\in[0,1]} \int_{x_0-r_j}^{x_0-\frac12r_j}(r_j^2-(x-x_0)^2)\,dx \geq c\frac{h_j^2}{j}r_j^2\\ \geq c\frac1j\Big(\frac{h_j}{k_j}\Big)^2\geq cj.
\end{multline}
where $c\in(0,1)$ is a small fixed quantity.\hfill{$\square$}

By using the method of proof above, it is clear that we may also directly show that $\n P(A)=+\infty$.

By virtue of Theorem \ref{thm.main}, we may deduce heuristically that the $A_{\infty}$ constant of $\omega_j^P$ must blow up as $j\ra\infty$. We now present a rigorous description of this fact.

\begin{proposition}[Degeneracy of the quantitative absolute continuity]\label{prop.deg} For any $q>1$, $|\m K_j|_{RH_q}$ goes to infinity as $j\ra\infty$. Moreover, $|\m K_j|_{A_{\infty}}$ goes to infinity as $j\ra\infty$.
\end{proposition}

\noindent\emph{Proof.} We show the first statement. Suppose otherwise, so that there exists $q>1$ and a constant $C_q$ such that for each $j\in\bb N$ and each surface ball $\Delta\subset\partial\Omega$, the estimate (\ref{eq.rhcond}) holds with $w\equiv\m K_j$ and $C\equiv C_q$. In particular, by setting $\Delta=[0,1]\times\{0\}$, we have that (\ref{eq.rhcond}) reduces to
\[
\Vert\m K_j\Vert_{L^q[0,1]}\leq C_q\Vert\m K_j\Vert_{L^1[0,1]},\qquad j\in\bb N.
\]
Using (\ref{eq.likerieszk}), the above estimate implies that
\begin{equation}\label{eq.compact}
\Vert\m R_j\Vert_{L^q[0,1]}\lesssim_q\Vert\m R_j\Vert_{L^1[0,1]},\qquad j\in\bb N.
\end{equation}
In fact, we have that (see \cite{cfa} page 233) 
\begin{equation}\label{eq.l11}
\Vert\m R_j\Vert_{L^1[0,1]}=1,\qquad j\in\bb N.
\end{equation}
Therefore, from (\ref{eq.compact}) we deduce that 
\begin{equation}\label{eq.bddq}
\Vert\m R_j\Vert_{L^q[0,1]}\lesssim_q1,
\end{equation}
which in particular implies that the family $\{\m R_j\}$ is uniformly integrable on $[0,1]$. Recall that (\cite{zyg}, Chapter V, Section 7, Theorem 7.7) $\m R_j\ra0$ pointwise a.e. on $[0,1]$. Then by (\ref{eq.bddq}) and the de la Vall\'ee Poussin criterion for equiintegrability (see \cite{bogachev} Volume I, Theorem 4.5.9), we must conclude by the Vitali Convergence Theorem that $\m R_j\ra0$ in $L^1[0,1]$ as $j\ra\infty$, but this stands in direct contradiction to (\ref{eq.l11}). The first desired statement follows.

By the same technique as above, we can verify that $|\n K_j|_{RH_{L\log L}}\ra\infty$ as $j\ra\infty$, and this must imply (quantitatively \cite{br}, actually) that $|\n K_j|_{A_{\infty}}\ra\infty$. \hfill{$\square$}

Next, let us tweak some parameters in the Modica-Mortola construction to obtain 

\begin{proposition}[Degeneracy of the Kenig-Pipher condition while $A_{\infty}$ is preserved]\label{prop.kenigbreaksdown} There exists a sequence $\{A_j\}$ of diagonal elliptic matrices on $\Omega$ with smooth, bounded, real coefficients in $\Omega$ and uniformly continuous on $\overline\Omega$, and there exists a diagonal elliptic matrix $A$ on $\Omega$ with smooth, bounded, real coefficients in $\Omega$ and uniformly continuous on $\overline\Omega$, such that $A_j\ra A$ pointwise uniformly on $\overline\Omega$, $\omega_A$ is absolutely continuous with respect to the surface measure $\sigma$, and
\begin{equation}\label{eq.result2}
\n P(A_j)\gtrsim j,\qquad\n P(A)=+\infty.
\end{equation}
\end{proposition}

\noindent\emph{Proof.} Consider the following modifications: First, we use the formula
\[
\phi_j=1+\frac1{4\pi j}\cos(2\pi h_jx),\qquad j\in\bb N
\]
for $\phi_j$ instead of the formula (\ref{eq.phij}) in Section \ref{sec.mm}. Second, we ask that the lacunary sequences $\{h_j\}$, $\{k_j\}$ satisfy the additional stronger estimate
\[
h_j\geq j^3 k_j,\qquad j\in\bb N.
\]
See the footnote to (\ref{eq.lac}).

With these changes in mind and following the argument for the proof of Theorem \ref{thm.main}, we still conclude as before that the measures $\omega_j$ converge weakly to a measure $\omega$, and that the corresponding Riesz products $\m R_j(x)=\prod_{i=1}^j\phi_i$ form the partial sums of a Fourier-Stieltjes series for a non-decreasing continuous function $F$ on $[-1,1]$. Moreover, we may mimic the proof of Proposition \ref{prop.break} and easily deduce (\ref{eq.result2}) accordingly.

On the other hand, in this situation, the amplitude coefficients $\frac1{2j}$ of the Riesz Products   are such that the sum of their squares is finite. According to \cite{zyg} Chapter V, Section 6, Lemma 6.5, it follows that $\{\m R_j\}$ is a uniformly bounded sequence in $L^2[0,1]$ which converges pointwise a.e. on $[0,1]$ to $F'$. Per the Vitali Convergence Theorem, we must conclude that $\m R_j\ra F'$ strongly in $L^1[0,1]$, whence the Fundamental Theorem of Calculus applies. Consequently, $dF=F'dx$ on $[0,1]$. Since $\m K_j\approx\m R_j$, it finally follows that $\omega\ll dx$, and we do remark that $\omega\in A_{\infty}$ holds via the Vitali Convergence Theorem and the Reverse H\"older inequality. \hfill{$\square$}

\begin{remark}\label{rm.top} Let us reframe the above results as follows. Define $\bb M$ as the Fr\'echet space of $n\times n$ matrix functions in $C^0(\overline{\Omega})\cap C^{\infty}(\Omega)$ with the usual topology. Designate $\operatorname{KP}\subset\bb M$ as the subset of such matrices which satisfy the Kenig-Pipher condition, and \footnote{As is well-known, such matrices are regular for the Dirichlet problem; this is the impetus for our notation.} $\operatorname{D}\subset\bb M$ as the subset of such matrices for whom the associated elliptic measure lies in $A_{\infty}$. Let $\partial=\partial_{\bb M}$ be the boundary operator on $\bb M$.
	
In this setting, note that the theorem of Kenig-Pipher \cite{kp3} is  the statement that $\operatorname{KP}\subset\operatorname{D}$. Hence $\partial\operatorname{KP}\subset\overline{\operatorname{D}}$. What we have done in the previous propositions is to parse the relationship between $\partial\operatorname{KP}$ and $\overline{\operatorname{D}}$ more delicately. Indeed, observe that for $A\in\partial\operatorname{KP}$, we must necessarily have $\n P(A)=+\infty$, and if $A'\in\partial\operatorname{D}$, then it must be the case that $\omega_{A'}$ is not $A_{\infty}$. Proposition \ref{prop.break} gives that 
\[
\partial\operatorname{KP}\cap\partial\operatorname{D}\neq\varnothing,
\]
which is not too surprising in light of the \cite{cfk} and \cite{mm} examples. On the other hand, Proposition \ref{prop.kenigbreaksdown} yields that
\[
\partial\operatorname{KP}\cap \operatorname{D}\neq\varnothing,
\]
which prohibits any quantifiable ``equivalence'' between the Kenig-Pipher condition and the $A_{\infty}$ property.

\end{remark}

\bibliography{refs}{} 

\begin{thebibliography}{CFMS81}

\bibitem[Bog07]{bogachev}
V.~I. Bogachev.
\newblock {\em Measure theory. {V}ol. {I}, {II}}.
\newblock Springer-Verlag, Berlin, 2007.
\newblock \href {https://doi.org/10.1007/978-3-540-34514-5}
  {\path{doi:10.1007/978-3-540-34514-5}}.

\bibitem[BR14]{br}
O.~Beznosova and A.~Reznikov.
\newblock Sharp estimates involving {$A_\infty$} and {$L\log L$} constants, and
  their applications to {PDE}.
\newblock {\em Algebra i Analiz}, 26(1):40--67, 2014.
\newblock URL:
  \url{https://doi-org.ezp1.lib.umn.edu/10.1090/s1061-0022-2014-01329-5}, \href
  {https://doi.org/10.1090/s1061-0022-2014-01329-5}
  {\path{doi:10.1090/s1061-0022-2014-01329-5}}.

\bibitem[CFK81]{cfk}
Luis~A. Caffarelli, Eugene~B. Fabes, and Carlos~E. Kenig.
\newblock Completely singular elliptic-harmonic measures.
\newblock {\em Indiana Univ. Math. J.}, 30(6):917--924, 1981.
\newblock \href {https://doi.org/10.1512/iumj.1981.30.30067}
  {\path{doi:10.1512/iumj.1981.30.30067}}.

\bibitem[CFMS81]{cfms}
L.~Caffarelli, E.~Fabes, S.~Mortola, and S.~Salsa.
\newblock Boundary behavior of nonnegative solutions of elliptic operators in
  divergence form.
\newblock {\em Indiana Univ. Math. J.}, 30(4):621--640, 1981.
\newblock \href {https://doi.org/10.1512/iumj.1981.30.30049}
  {\path{doi:10.1512/iumj.1981.30.30049}}.

\bibitem[Dah77]{dah1}
Bj\"orn E.~J. Dahlberg.
\newblock Estimates of harmonic measure.
\newblock {\em Arch. Rational Mech. Anal.}, 65(3):275--288, 1977.
\newblock \href {https://doi.org/10.1007/BF00280445}
  {\path{doi:10.1007/BF00280445}}.

\bibitem[Gra14]{cfa}
Loukas Grafakos.
\newblock {\em Classical {F}ourier analysis}, volume 249 of {\em Graduate Texts
  in Mathematics}.
\newblock Springer, New York, third edition, 2014.
\newblock \href {https://doi.org/10.1007/978-1-4939-1194-3}
  {\path{doi:10.1007/978-1-4939-1194-3}}.

\bibitem[KP01]{kp3}
Carlos~E. Kenig and Jill Pipher.
\newblock The {D}irichlet problem for elliptic equations with drift terms.
\newblock {\em Publ. Mat.}, 45(1):199--217, 2001.
\newblock \href {https://doi.org/10.5565/PUBLMAT_45101_09}
  {\path{doi:10.5565/PUBLMAT_45101_09}}.

\bibitem[MM81]{mm}
Luciano Modica and Stefano Mortola.
\newblock Construction of a singular elliptic-harmonic measure.
\newblock {\em Manuscripta Math.}, 33(1):81--98, 1980/81.
\newblock \href {https://doi.org/10.1007/BF01298340}
  {\path{doi:10.1007/BF01298340}}.

\bibitem[Ste93]{ste}
Elias~M. Stein.
\newblock {\em Harmonic analysis: real-variable methods, orthogonality, and
  oscillatory integrals}, volume~43 of {\em Princeton Mathematical Series}.
\newblock Princeton University Press, Princeton, NJ, 1993.
\newblock With the assistance of Timothy S. Murphy, Monographs in Harmonic
  Analysis, III.

\bibitem[Zyg02]{zyg}
A.~Zygmund.
\newblock {\em Trigonometric series. {V}ol. {I}, {II}}.
\newblock Cambridge Mathematical Library. Cambridge University Press,
  Cambridge, third edition, 2002.
\newblock With a foreword by Robert A. Fefferman.

\end{thebibliography}
\bibliographystyle{alphaurl} 

\end{document}